\documentclass[12pt,letterpaper]{amsart}
\usepackage{amssymb,amsmath,amscd,eucal,epsfig}
\def\<{\langle}                     %added
\def\>{\rangle}                     %added
\newcommand{\ben}{\begin{enumerate}}
\newcommand{\een}{\end{enumerate}}

\theoremstyle{plain}

\newtheorem{proposition}{Proposition}[section]

\newtheorem{remark}{Remark}[section]

\theoremstyle{definition}
\newtheorem{definition}{Definition}[section]

\newtheorem{example}{Example}[section]

\numberwithin{equation}{section}
\begin{document}

\begin{center}
{\textbf{{An Application of Complex Fuzzy Soft Matrices in Signal Processing}}}\ \\ \ \\

{Olayemi R. Oladokun}\\ olayemirita@gmail.com\\ and \\ Taiwo O. Sangodapo$^*$ \\ to.sangodapo@ui.edu.ng\\
Department of Mathematics,\\ University of Ibadan, Ibadan, Nigeria\\Corresponding author $*$ 
\end{center}
\ \\ \  \\

\begin{center} {\textbf{Abstract}} \end{center}
In this paper, we study the concept of complex fuzzy soft matrices. The application of complex fuzzy soft matrices in signals and systems via the cross product of complex fuzzy soft matrices and Fourier transform was carried out. In this application, an algorithm for the identification of a reference signal out of large interest signals detected by a digital receiver was presented. It was recorded that, the Fourier transform is better because it gave a higher optimal value and as a result, there was a better reference signal $R.$
\; \\ \; \\ \; \\

\ \\


\section{Introduction}
Most of our traditional tools for formal modelling, reasoning, and computing are crisp, deterministic, and precise in character. By crisp we mean dichotomous, that is, yes-or-no-type rather than more-or-less type. In conventional dual logic, for instance, a statement can be true or false-and nothing in between. In set theory, an element can either belong to a set or not; and in optimization, a solution is either feasible or not. Precision assumes that the parameters of a model represent exactly either our perception of the phenomenon modelled or the features of the real system that has been modelled. Certainty eventually indicates that we assume the structures and parameters of the model to be definitely known, and that there are no doubts about their values or their occurrence. \\If the model under consideration is a formal model Zimmermann in 1980 then the model assumptions are in a sense arbitrary, that is, the model builder can freely decide which model characteristics he chooses. If, however, the model or theory asserts factuality Zimmermann in 1980 then the modelling language has to be suited to model the characteristics of the situation under study appropriately.

Fuzzy set has been introduced by Zadeh in 1965 as a generalisation of classical set. This theory is a suitable tool for modeling since the classical models lack clarity while handling problems in different fields of science like artificial intelligence, decision theory, and many others. 
In other words, this concept of fuzzy set is used to address those uncertain problems which arise in models representing real life phenomenon. This theory has been studied and applied by many researchers, see \cite{ta1, ta2}.

In 2013, a generalisation of fuzzy set was introduced by Cuong and Kreinovich \cite{ck} to deal with uncertainty and imprecision in many decision making problems. Their work has been studied, applied and extended, see \cite{cfs, dg, t1, t2, t3, to}.   

In 1999, Molodsov \cite{m} extended the notion of fuzzy set to soft set theory. In 2009, Ali et al \cite{a9}, studied the soft set and in 2010, Cagman and Enginoglu \cite{a15} also studied soft set theory and applied it to decision making problem. This theory was extended to fuzzy soft matrix by Cagman and Enginoglu \cite{a4, a16, a17} in 2012 and 2020, respectively. In 2021, Sangodapo, Onasanya and Hoskova-Mayerova \cite{soh} applied fuzzy soft matrix to medical diagnosis using a revised method. Sangodapo and Feng \cite{sf} in 2022, applied fuzzy matrix composition to medical diagnosis via an improved method.  

Ramot et al. in 2012, pioneered the notion of complex fuzzy set theory which is the generalisation of fuzzy set.  whose range is not restricted to a closed interval $[0, 1]$ but expanded to a unit circle in complex plane. This theory is different from the theory of fuzzy set because it has a phase term and amplitude term. In 2015. Ahmed and Mahanta \cite{a5} introduced the concept of complex fuzzy soft matrix and studied their properties and operations. Haque and Kar \cite{a20} and Madad et al, 2021 and 2020, respectively studied complex fuzzy soft matrices and applied it to decision making problems.

In this paper, the application of complex fuzzy soft sets in signals and systems via the cross product of complex fuzzy soft matrices and Fourier transform was carried out. In this application, an algorithm for the identification of a reference signal out of large interest signals detected by a digital receiver was presented. It was recorded that, the Fourier transform is better because it gave a higher optimal value and as a result, there was a better reference signal $R.$

\ \\

\section{Preliminaries}
This section gives the basic definitions and existing results relating to complex fuzzy soft sets and matrices.

a\begin{definition} \cite{m}
A pair $(F,E)$ is called a soft set $(over\; U)$ if and only if $F$ is a mapping of $E$ into the set of all subsets $P(U)$ of the set $U$.
\end{definition}

In other words, the soft set is a parametrized family of subsets of the set $U$. Every set $F(c),\;c\in E,$ from this family may be considered as the set of $c$-elements of the soft set $(F, E),$ or as the set of $c$-approximate elements of the soft set.

\begin{example}
A soft set $(F, E)$ describes the attractiveness of the houses which Mr. $X$ is going to buy.\\
$U$ - is the set of houses under consideration.\\
$E$ - is the set of parameters. Each parameter is a word or a sentence.\\
$E = \{expensive;\; beautiful;\; wooden; \;cheap;\; in\; the \;green\; surroundings;\\ modern;\; in\; good\;
repair;\; in \;bad \;repair\}.$
\end{example}

\begin{definition}\cite{ma}
Let $U$ be an initial universe set, $E$ be the set of parameters and $A_i\subset E$. Then a pair $(\Phi_{\mu},A)$ is called fuzzy soft set over $U$, where $\Phi_{\mu}$ is a mapping given by $$\Phi_{\mu}:A_i\to P(\mu_U),$$ where $P(\mu_U)$ denote the set of complex fuzzy subsets of $U$.\\
\end{definition}

\begin{definition}\cite{ma}
For two fuzzy soft sets $(F,A)$ and $(G,B)$ over a common universe $U,$ $(F,A)$ is a fuzzy soft subset of
$(G,B)$ if\\
$(i)\;\; A\subset B,$ and\\
$(ii)\;\;\forall c\in A,\; F(c)$ is a fuzzy subset of $G(c).$ We write $(F,A)\tilde\subset(G,B).$\\
$(F,A)$ is said to be a fuzzy soft super set of $(G,B),$ if $(G,B)$ is a fuzzy soft subset of $(F,A).$ denoted by
$(F,A)\tilde\supset(G,B).$
\end{definition}

\begin{example}
Consider two fuzzy-soft-sets $(F,A)$ and $(G,B)$ over the same universal set $U = \{h_1, h_2, h_3, h_4, h_5\}.$ Here
$U$ represents the set of houses, $A = \{blackish, reddish, green\}$ and $B = \{blackish, reddish, green, large\},$ and for $F$ and $G$ we have\\
$F(blackish) = \{h_1/.4, h_2/.6, h_3/.5, h_4/.8, h_5/1\},$\\
$F(reddish) = \{h_1/1, h_2/.5, h_3/.5, h_4/1, h_5/.7\},$\\
$F(green) = \{h_1/.5, h_2/.6, h_3/.8, h_4/.8, h_5/.7\},$\\
$G(blackish) = \{h_1/.4, h_2/.7, h_3/.6, h_4/.9, h_5/1\},$\\
$G(reddish) = \{h_1/1, h_2/.6, h_3/.5, h_4/1, h_5/8\},$\\
$G(green) = \{h_1/.6, h_2/.6, h_3/.9, h_4/.8, h_5/7\},$\\
$G(large) = \{h_1/.4, h_2/.6, h_3/.5, h_4/.8, h_5/1\}.$\\
Clearly, $(F,A)\tilde\subset(G,B).$
\end{example}
\ \\ \ \\

\section{Fuzzy Soft Matrices}
In this section, we define fuzzy soft matrices which are representations of the fuzzy soft sets.\\
This style of representation is useful for storing a soft set in a computer memory. The operations can be presented by the matrices which are very useful and  convenient for application.\\
The set of all fuzzy sets over $U$ will be denoted by $F(U).$ $\tau_A$, $\tau_B$, $\tau_C,\cdots,$ etc and $\phi_A$, $\phi_B$, $\phi_C,\cdots,$ etc will be used for fuzzy soft sets and their fuzzy approximate functions, respectively.

\begin{definition}\cite{a16}
Let $U$ be an initial universe, $E$ be the set of all parameters, $A\subseteq E$ and 
$\phi_A(x)$ be a fuzzy set over $U$ for all $x\in E$. Then, a fuzzy soft set $\tau_A$ over $U$ is a set defined by a function 
$\phi_A$ representing a mapping
$$\phi_A:E\to F(U)$$ such that $$\phi_A(x)=\emptyset\;\;if\;x\notin A$$
Here,  $A$ is called fuzzy approximate function of the fuzzy soft set $\tau_A$, the value  $\phi_A(x)$
is a fuzzy set called $x$-element of the fuzzy soft set for all $x\in E,$ and is the null fuzzy
set. Thus, a fuzzy soft set $\tau_A$ over $U$ can be represented by the set of ordered pairs
$$\tau_A=\{(x,\phi_A(x)): x\in E,\phi_A(x)\in F(U)\}.$$
\end{definition}
the sets of all fuzzy soft sets over $U$ will be denoted by $FS(U).$

\begin{example}
Assume that $U=\{u_1,u_2,u_3,u_4,u_5\}$ is a universal set and $E=\{x_1,x_2,x_3,x_4\}$ is a set of all parameters.\\ If $A=\{x_2,x_3,x_4\}$, $\phi_A(x_2)=\{0.5/u_2, \;0.8/u_4\},$, $\phi_A(x_3)=\emptyset$ and $\phi_A(x_4)=U,$ then the fuzzy soft set $\tau_A$ is written by $$\tau_A=\{(x_2,\{0.5/u_2,\;0.8/u_4\}),(x_4,U)\}.$$
\end{example}

\begin{definition}
Let $\tau_A\in FS(U).$ Then, a fuzzy relation form of $\tau_A$ is defined by $$R_A=\{(\mu_{R_A}(u,x)/(u,x):(u,x)\in U\times E\},$$
where the membership function of $\mu_{R_A}$ is written by
$$\mu_{R_A}=:U\times E\to [0,1],\;\;\;\mu_{R_A}(u,x)=\mu_{\phi_A(x)}(u).$$
\end{definition}

If $U=\{u_1,u_2,\cdots,u_m\}$,$\;\;E=\{x_1,x_2,\cdots,x_n\}$ and $A\subseteq E$, then the $R_A$ can be presented by a table as in the following form
$$\begin{tabular}{c|ccccc}
$R_A$&$x_1$&$x_2$&$\cdots$ &$x_n$\\\hline
$u_1$&$\mu_{R_A}(u_1,x_1)$&$\mu_{R_A}(u_1,x_2)$&$\cdots$&$\mu_{R_A}(u_1,x_n)$\\
$u_2$&$\mu_{R_A}(u_2,x_1)$&$\mu_{R_A}(u_2,x_2)$&$\cdots$&$\mu_{R_A}(u_2,x_n)$\\
$\vdots$&$\cdots$&$\cdots$&$\ddots$&$\vdots$\\
$u_m$&$\mu_{R_A}(u_m,x_1)$&$\mu_{R_A}(u_m,x_2)$&$\cdots$&$\mu_{R_A}(u_m,x_n)$
\end{tabular}$$
If $a_{ij}=\mu_{R_A}(u_i,x_j),$ we can define a matrix 
$$[a_{ij}]_{m\times n}=\left[\begin{array}{cccc}
x_{11}&x_{12}&\cdots& x_{1n}\\
x_{21}&x_{22}&\cdots&x_{2n}\\
\vdots&\vdots&\ddots&\vdots\\
x_{m1}&x_{m2}&\cdots&x_{mn}
\end{array}\right]$$
which is called an $m\times n$ fuzzy soft matrix of the fuzzy soft set $\tau_A$ over $U.$

\begin{remark}
According to this definition, a fuzzy soft set $\tau_A$ is uniquely characterized by the matrix $[a_{ij}]_{m\times n}.$ It means that a fuzzy soft set $\tau_A$ is formally equal to its soft matrix $[a_{ij}]_{m\times n}.$ Therefore, we shall identify any fuzzy soft set with its fuzzy soft matrix and use these two concepts interchangeable.
\end{remark}

The set of all $m\times n$ fuzzy soft matrices over $U$ will be denoted by $FSM_{m\times n}.$ From now
on we shall delete the subscript $m\times n$ of $[a_{ij}]_{m\times n}$, we use $[a_{ij}]$ instead of$[a_{ij}]_{m\times n}$,
since $[a_{ij}]_{m\times n}\in FSM_{m\times n}$ means that $[a_{ij}]$ is an $m\times n$ fuzzy soft matrix for $i = 1, 2,\cdots,m$ and $j = 1, 2,\cdots,n.$

\begin{example}
Let us consider the table in the previous definition. Then, the relation form of $\tau_A$ is written by
$$R_A=\{0.5/(u_2,x_2),\;0.8/(u_4,x_2),\;1/(u_1,x_4),\;1/(u_2,x_4),\;1/(u_3,x_4),\;1/(u_4,x_4),\;1/(u_5,x_4)\}$$ Hence, the fuzzy soft matrix $[a_{ij}]$ is written by 
$$[a_{ij}]=\left[\begin{array}{cccc}
0 & 0 & 0 &1\\
0 & 0.5 & 0 & 1\\
0 & 0 & 0 & 1\\
0 & 0.8 & 0 & 1\\
0 & 0 & 0 & 1
\end{array}\right]$$
\end{example}

\begin{definition}\cite{a16}
Let $[a_{ij}],[b_{ik}]\in FSM_{m\times n}.$ Then, And product of $[a_{ij}]$ and $[b_{ik}]$ is defined by $$\land:FSM_{m\times n}\times FSM_{m\times n}\to FSM_{m\times n^2},\;\;\;[a_{ij}]\land[b_{ik}]=[c_{ip}]$$ where $c_{ip}=min\{a_{ij},b_{ik}\}$ such that $p=n(j-1)+k.$
\end{definition}

\begin{definition} \cite{a16}
Let $[a_{ij}],[b_{ik}]\in FSM_{m\times n}.$ Then, Or product of $[a_{ij}]$ and $[b_{ik}]$ is defined by $$\lor:FSM_{m\times n}\times FSM_{m\times n}\to FSM_{m\times n^2},\;\;\;[a_{ij}]\lor[b_{ik}]=[c_{ip}]$$ where $c_{ip}=max\{a_{ij},b_{ik}\}$ such that $p=n(j-1)+k.$
\end{definition}

\begin{definition} \cite{a16}
Let $[a_{ij}],[b_{ik}]\in FSM_{m\times n}.$ Then, And Not product of $[a_{ij}]$ and $[b_{ik}]$ is defined by $$\overline{\land}:FSM_{m\times n}\times FSM_{m\times n}\to FSM_{m\times n^2},\;\;\;[a_{ij}]\overline{\land}[b_{ik}]=[c_{ip}]$$ where $c_{ip}=min\{a_{ij},1-b_{ik}\}$ such that $p=n(j-1)+k.$
\end{definition}

\begin{definition} \cite{a16}
Let $[a_{ij}],[b_{ik}]\in FSM_{m\times n}.$ Then, Or Not product of $[a_{ij}]$ and $[b_{ik}]$ is defined by $$\overline{\lor}:FSM_{m\times n}\times FSM_{m\times n}\to FSM_{m\times n^2},\;\;\;[a_{ij}]\overline{\lor}[b_{ik}]=[c_{ip}]$$ where $c_{ip}=max\{a_{ij},1-b_{ik}\}$ such that $p=n(j-1)+k.$
\end{definition}
\ \\ \ \\

\section{Complex Fuzzy Matrices}
\begin{definition} \cite{a23}
A complex fuzzy matrix is a matrix which has its elements from the unit disc in the complex plane. A complex fuzzy matrix of order $m \times n$ is defined as, $$A=[\langle a_{ij},a_{ij\mu}\rangle]_{m\times n}a_{ij\mu}=r_{ij}e^{i\omega_{ij}}$$ where $r_{ij}\in[0,1]\;and\;\omega_{ij}\in [0,2\pi)$ \\
\end{definition}

\textbf{Some Properties of Complex Fuzzy Matrices.}\\
1.\;\;Two complex fuzzy matrices are conformable for addition if the matrices are of the same order. If $A=[a_{ij\mu}]_{m\times n}$ and $B=[b_{ij\theta}]_{m\times n}$ then
$$A+B=[c_{ij\nu}]_{m\times n};\;\;c_{ij\nu}=max\{a_{ij\mu},\;b_{ij\theta}.$$

\begin{example} \label{a}
Let $$A=\left[\begin{array}{ccc} 
0.6e^{i{\pi\over2}}&0.4e^{i{\pi\over2}}&0.1e^{i0}\\
0.0e^{i0}&0.1e^{i\pi}&0.3e^{i{\pi\over4}}\\
1.0e^{i0}&0.2e^{i{\pi\over3}}&0e^{i0}
\end{array}\right]\;\;and \;\;
B=\left[\begin{array}{ccc}
0.1e^{i{\pi\over6}}&0.2e^{i\pi}&0.5e^{i0}\\
0.8e^{i0}&0.4e^{i{\pi\over4}}&0.7e^{i\pi}\\
0.3e^{i{\pi\over4}}&0.0e^{i0}&1.0e^{i\pi}
\end{array}\right]$$

then
$$A+B=\left[\begin{array}{ccc} 
0.6e^{i{\pi\over2}}&0.4e^{i{\pi\over2}}&0.5e^{i0}\\
0.8e^{i0}&0.4e^{i{\pi\over4}}&0.7e^{i\pi}\\
1.0e^{i{\pi\over4}}&0.2e^{i{\pi\over3}}&1.0e^{1\pi}
\end{array}\right]$$  
\end{example}

2.\;\;The product of two complex fuzzy matrices under usual matrix multiplication is not a complex fuzzy matrix. So here we intend to define a compactible operation analogous to product that the product again happens to be a complex fuzzy matrix. To define the product, the number of columns of the first matrix must be equal to the number of rows of the second matrix. i.e.; If $A=[a_{ij\mu}]_{m\times n}$ and $B=[b_{ij\theta}]_{n\times p}$

$$AB=[c_{ij\nu}]_{m\times p};$$$\;\;c_{ij\nu}=max\{min\{a_{i1\mu},\;b_{1j\theta}\}\},max\{min\{a_{i2\mu},\;\;b_{2j\theta}\}\},\; max\{min\{a_{i3\mu},\;b_{3j\theta}\}\},\\\cdots,max\{min\{a_{nj\mu},\;b_{in\theta}\}\},\;\;1\leq i\leq m,\;1\leq j\leq p\;\;\;$

\begin{example}
From Example \ref{a}, $$AB=\left[\begin{array}{ccc} 
	0.4e^{i{\pi\over2}}&0.4e^{i{\pi\over2}}&0.5e^{i0}\\
	0.3e^{i0}&0.1e^{i{\pi\over4}}&0.3e^{i\pi}\\
	0.2e^{i{\pi\over4}}&0.2e^{i{\pi\over3}}&0.5e^{1\pi}
\end{array}\right]$$ 
\end{example}

3\;\;Let $A=[a_{ij\mu}]_{n\times n}$ be a square complex fuzzy matrix. Then, Trace(A)=$max\{a_{ii\mu}\}$.

\begin{example}
From Example \ref{a}, Trace(A) = $0.6$
\end{example}

4.\;\;The conjugate transpose of a complex fuzzy matrix $A$ is defined as the transpose of the conjugate of the matrix.
That is
$$\overline{A}^T=(A)^{\ast}$$

\begin{definition} \cite{a23}
Let $U = \{c_1, c_2, c_3,\cdots , c_n\}$ be the universal set and $E$ be the set of parameters given by $E = \{e_1, e_2, e_3, \cdots , e_m\}$. If $A\subset E,$ then $(\Phi_{\mu},A)$ is a fuzzy soft matrix over $U,$ where $\Phi_{mu}$ is a mapping given by
$$\Phi_{\mu}: A\to P(\mu_U)\;and\;
\nu_{\mathcal{E}}: A\to C^U$$ where $C^U$ denotes the set of all fuzzy subsets and $P(C^U)$ denotes the
power set of $C^U.$ Then, complex fuzzy soft set $(\Phi_{\mu},A)$ can be expressed in matrix form as
$$[\overline{A}_{m\times n}] = [|a_{ij}|]_{m\times n},$$ for $i = 1, 2, \cdots ,m$ and $j = 1, 2, 3,\cdots , n,$ where
$$|a_{ij}|=\begin{cases} |\nu(e_i)|_j,\;\;\;\;\;if\;e_i\in A\\0,\;\;\;\;\;\;\;\;otherwise \end{cases}$$
$|\nu(e_i)|_j\;\;\;(\{\nu(e_i)\}_j$ is a complex fuzzy set ) represents the element of $\overline{A}$ corresponds
to element $c_j$ of $U,$ for $j = 1, 2, 3,\cdots , n,$ where $|\nu(e_i)|_j = \alpha_{RC}$ such that $\alpha\in[0, 1]$ and
$R = 1, 2, 3, \cdots ,m$ and $C = 1, 2, 3, \cdots , n.$	
\end{definition}

\begin{example} \label{2}
Suppose that there are three houses under consideration, namely the universes $U = \{h_1, h_2, h_3\},$ and the parameter set $E = \{e_1, e_2, e_3\},$ where ei stands for “near to downtown”, “green”, and “cheap ” respectively. Consider the mapping $\Phi_{\mu}$ from parameter set $A = \{e_1, e_2\}\subset E$ to the set of all fuzzy subsets of power set $U.$ Consider a
fuzzy soft set $(\Phi_{\mu},A)$ which describes the “attractiveness of houses” that is considering for purchase. Then, the fuzzy soft set $(\Phi_{\mu},A)$ is given as.
$$[\Phi_{\mu},A]=[|\nu(e_i)|_j]_{m\times n},$$ where 
$\Phi_{\mu}(e_1)=\{(h_1, |0.1e^{i{\pi\over2}} |), (h_2, |0.2e^{i{2\pi}}|, (h_3, |0.3e^{i{\pi\over4}}|)\},$\\
$\Phi_{\mu}(e_2)=\{(h_1, |0.4e^{i{\pi\over6}} |), (h_2, |0.5e^{i(0.5)}|, (h_3, |0.6e^{i(0)}|)\}$,\\
where $[note\; that \;e^{ix}=cos(x)+isin(x)],$\\
(1)\;\;\;$0.1e^{i{\pi\over2}}=0.1(cos{\pi\over2}+isin{\pi\over2})$,\\
$\;\;\;\;\;\;\;=0.1(cos{90} + isin{90}) = 0.1i$\\
But\\  
$|0.1i|^2 = |(0.1i)(0.1i)| = |0.01i^2| = |-0.01| = 0.01,$\\
Hence, $|0.1i| = \sqrt{0.01} = 0.1$\\ \ \\
(2)\;\;\;$0.2e^{i{\pi\over6}} = 0.2(cos{\pi\over6} + isin{\pi\over6})$,\\ $0.2e^{i{\pi\over6}} = 0.2(cos{30}+isin{30})$,
$ = 0.2(0.87 + i0.5) = 0.174 + i0.1$\\$|0.174 + i0.1| = \sqrt{0.03 + 0.01} = 0.2$\\ \ \\
(3)\;\;\;$0.3e^{i{2\pi}} = 0.3(cos{2\pi} + isin{2\pi}) = $
$= 0.3(0.99 + i0.05) = 0.297 + i0.015$\\$|0.297 + i0.015| = \sqrt{0.088+0.000} = 0.3$\\ \ \\
(4)\;\;\;$0.4e^{i{\pi\over4}} = 0.4(cos{cos45} + isin{sin45})$
$= 0.4(0.707 + i0.707) = 0.28 + i0.28$\\$|0.28 + i0.28| = \sqrt{0.0784 + 0.0784} = 0.16$\\ \ \\
(5)\;\;\;$0.5e^{i(0)} = 0.5(cos{0} + isin{0})$
$= 0.5$\\$|0.5| = \sqrt{0.25} = 0.5$\\ \ \\
(6)\;\;\;$0.6e^{i(0.5)}=0.6(cos{0.5}+isin{0.5}) = 0.6(0.99 + i0.009) = 0.594 + i0.005$\\$|0.594 + i0.005| = \sqrt{0.353} = 0.6$\\ \ \\
We would represent this complex fuzzy soft set in matrix form as:
$$\left(\begin{array}{ccc}
0.1 & 0.4 & 0\\
0.2 & 0.5 & 0\\
0.3 & 0.6 & 0
\end{array}\right).$$
\end{example}

\begin{definition} \cite{a23}
Let $[a_{ij}]$ be a complex fuzzy soft matrix. Then $[a^{{A}}_{ij}]$ is called complex zero soft matrix if $(a_{ij}, r_{ij}) = (0,0)$ for all i and j and denoted by h  $[a^{{A}}_{ij}]=[0].$\\
${A}_{m\times n}=[|a_{ij}|_{m\times n},$ for $i=1,2,\cdots,m$ and $j=1,2,\cdots,n,$ where
$$|a_{ij}|=\begin{cases} |\nu(e_i)|_j,\;\;\;\;\;if\;e_i\in A\\0,\;\;\;\;\;\;\;\;otherwise \end{cases}$$
$|\nu(e_i)|_j=\alpha_{RC}$ such that $\alpha\in [0,1]$ and $R=1,2,\cdots,m$ and $C=1,2,\cdots,n.$\\
$$\left(\begin{array}{ccc}
	0 & 0 & 0\\
	0 & 0 & 0\\
	0 & 0 & 0
\end{array}\right).$$
\end{definition}
\ \\
\begin{definition}\cite{a23}
Let $[a_{ij}]\in FSM_{m\times n}$. Then $[a_{ij}]$ is called\\
(1) a zero fuzzy soft matrix, denoted by $[0]$, if $a_{ij}=0$ for all $i$ and $j$.\\
(2) an $A$-universal fuzzy soft matrix, denoted by $[\overline{a}_{ij}]$, if $a_{ij}=1$ for all $i$ and $j\in I_A=\{j:x_j\in A\}$.\\
(3) a universal fuzzy soft matrix, denoted by $[1]$, if $a_{ij}=1$ for all $i$ and $j$.
\end{definition}

\begin{definition}\cite{a23}
Let $[{A}]$ and $[{B}]$ be complex fuzzy soft matrices. Then\\
(i)\;\;\;\; $[{A}]$ is a complex fuzzy soft submatrices of $[{B}]$, denoted by $[{A}]\subseteq [{B}]$,
if $|a_{ij}|\leq |b_{ij}|$ for all $|a_{ij}|\in[{A}]$ and $|b_{ij}|\in[{B}].$\\
(ii)\;\;\;\; $[{A}]$ is a proper complex fuzzy soft submatrices of $[{B}]$, denoted by $[{A}]\subset [{B}]$,
if $|a_{ij}|\leq |b_{ij}|$ for all $|a_{ij}|\in[{A}]$ and $|b_{ij}|\in[{B}], \;n,$ and for at least one term $|a_{ij}| < |b_{ij}|.$\\
(iii)\;\;\;\; $[{A}]$ is an equal complex fuzzy soft matrix of $[{B}_{m\times n}]$, denoted by $[{A}]= [{B}]$,
if $|a_{ij}|=|b_{ij}|$ for all $|a_{ij}|\in[{A}]$ and $|b_{ij}|\in[{B}].$
\end{definition}

\begin{definition}\cite{a23}
Let $[{A}]$ and $[{B}]$ be complex fuzzy soft matrices. Then the
complex fuzzy soft matrices $[{C}]$ are called\\
(i)\;\; Union of $[{A}]$ and $[{B}]$, denoted $[{A}]\cup[{B}]$. If $[{C}] = max\{|a_{ij}|, |b_{ij}|\}$ for all $|a_{ij}|\in[{A}]$ and $|b_{ij}| \in[{B}]$.\\
(ii)\;\; Intersection of $[{A}]$ and $[{B}]$, denoted $[{A}]\cap[{B}]$ if $[{C}] = min\{|a_{ij}|, |b_{ij}|\}$ for all $|a_{ij}|\in[{A}]$ and $|b_{ij}| \in[{B}]$.\\
(iii)\;\;Complement of $[{A}]$ denoted by $[{A}]^c$, if ${C}={1}- {A}$ for all $m$ and $n.$
\end{definition}
\ \\
\begin{example}
Following the same process in Example \ref{2}, we have these two complex fuzzy soft matrices $$A = \left(\begin{array}{cccc}
	0.1 & 0.4 & 0.0 & 0.2\\
	0.2 & 0.6 & 0.0 & 0.3\\
	0.7 & 0.2 & 0.0 & 0.5\\
	0.4 & 0.3 & 0.0 & 0.9
\end{array}\right)$$ and $$B = \left(\begin{array}{cccc}
0.2 & 0.4 & 0.0 & 0.1\\
0.3 & 0.7 & 0.0 & 0.2\\
0.8 & 0.8 & 0.0 & 0.4\\
0.5 & 0.3 & 0.0 & 0.7
\end{array}\right).$$ Then, $$A^c = \left(\begin{array}{cccc}
0.9 & 0.6 & 1.0 & 0.8\\
0.8 & 0.4 & 1.0 & 0.7\\
0.3 & 0.8 & 1.0 & 0.5\\
0.6 & 0.7 & 1.0 & 0.1
\end{array}\right),$$  $$A \cup B = \left(\begin{array}{cccc}
0.2 & 0.4 & 0.0 & 0.2\\
0.3 & 0.7 & 0.0 & 0.3\\
0.8 & 0.8 & 0.0 & 0.5\\
0.5 & 0.3 & 0.0 & 0.9
\end{array}\right)$$ and $$A \cap B = \left(\begin{array}{cccc}
0.1 & 0.4 & 0.0 & 0.1\\
0.2 & 0.6 & 0.0 & 0.2\\
0.7 & 0.2 & 0.0 & 0.4\\
0.4 & 0.3 & 0.0 & 0.7
\end{array}\right)$$ 
\end{example}
\ \\
\begin{proposition}\cite{a23}
	$[a_{ij}],[b_{ij}],[c_{ij}]\in CFSM_{m\times n}.$ Then
\begin{enumerate}
\item  $[a_{ij}]{\cap}[b_{ij}] = [b_{ij}]{\cap}[a_{ij}]$
\item $[a_{ij}]{\cup}[b_{ij}] = [b_{ij}]{\cup}[a_{ij}]$
\item $([a_{ij}]{\cap}[b_{ij}]){\cap}[c_{ij}] = [a_{ij}]{\cap}([b_{ij}]{\cap}[c_{ij}])$
\item $([a_{ij}]{\cup}[b_{ij}]){\cup}[c_{ij}] = [a_{ij}]{\cup}([b_{ij}]{\cup}[c_{ij}])$
\item $[a_{ij}]{\cap}([b_{ij}]{\cup}[c_{ij}]) = ([a_{ij}]{\cap}[b_{ij}]){\cup}([a_{ij}]{\cap}[c_{ij}])$
\item $[a_{ij}]{\cup}([b_{ij}]{\cap}[c_{i}j]) = ([a_{ij}]{\cup}[b_{ij}]){\cap}([a_{ij}]{\cup}[c_{ij}])$
\item $([a_{ij}{\cup}][b_{ij}])^{\circ} =[a_{ij}]^{\circ}{\cap}[b_{ij}]^{\circ}$
\item $([a_{ij}{\cap}][b_{ij}])^{\circ} =[a_{ij}]^{\circ}{\cup}[b_{ij}]^{\circ}$
\end{enumerate}
\end{proposition}
\ \\ \ \\

\section{Decision Making Algorithm in Signal Processing using Complex Fuzzy Soft Matrices}
We defined a real life application of newly defined complex fuzzy soft matrix and showed how our theoretical results have real life applications. The complex fuzzy soft matrix explains how to get better and clear signal for identification with given reference signal.

\begin{definition}
Let $U = \{u_1, u_2, u_3,\cdots u_n\}$ be initial universal set and $Mm(C_{ki}) =(D_{i1}).$ Then subset of $U$ can be obtained by using $[D_{i1}] =\overline{opt}[(\overline{d}_{i1},d_{i1})],$ where
$\overline{opt}[(\overline{d}_{i1})(U)=\{\overline{d}_{i1}/u_i:\;u_i\in U\}$ and ${opt}[({d}_{i1})(U)=\{{d}_{i1}/u_i:\;u_i\in U\},$
$\{\overline{opt}(\overline{d}_{i1}),{opt}({d}_{i1}U)\}$ is called an optimum fuzzy set on $U.$
\end{definition}

{\bf Step 1.}
If a receiver gets various signals $S_{1(g)}, S_{2(g)}, S_{3(g)},\cdots, S_{m(g)}$ from any source. Each signal is sampled N times by the receiver. Then $S_{i(g)}$ ($i$ varies from $1$ to $m$) signals can be recognized with respect to $R,$ where $R$ is given known signal. Assume that both $S_{i(g)}$ and $R$ are considered as $n$ times.\\
Assume that $S_{j(t)}$ is $i^{th}$ signal, where $1\leq t\leq m.$ Then the absolute value of each $S_{j(t)}$ in terms of discrete complex fuzzy transform given as follows: where $g_{j,n}$ is the complex Fourier coefficients of signals $S_j$ , and $S_{j(t)}$ varies from $1$ to $m.$ The above expression with alternate complex Fourier coefficients can be written as:
where $g_{j,s}=u_{j,s}\cdot e^{i\beta_{j,s}}$ such that $u_{j,s}$ and $\beta_{j,s}$ are real valued and $u_{j,s}\geq1$ for all $s\;(1 \leq s\leq m).$\\
{\bf Step 2.}
Expressed in matrix form as $\overline{A}_{m\times n} = [|S_j(t)|]_{N\times m},$ that is, in the matrix take all the signals in columns and each column contains $N$ samples of every signal, so we get

$$A=\left[\begin{array}{cccc}
|S_1(1)|&|S_2(1)|&\cdots&|S_m(1)|\\
|S_1(2)|&|S_2(2)|&\cdots&|S_m(2)|\\
\vdots&\vdots&\ddots&\vdots\\
|S_1(N)|&|S_2(N)|&\cdots&|S_m(N)|
\end{array}\right].$$

{\bf Step 3.} In same way make another matrix by the signals $S'_j(t)\; (1\leq t\leq m\; and \;1 \leq j\leq N).$
$$B=\left[\begin{array}{cccc}
|S'_1(1)|&|S'_2(1)|&\cdots&|S'_m(1)|\\
|S'_1(2)|&|S'_2(2)|&\cdots&|S'_m(2)|\\
\vdots&\vdots&\ddots&\vdots\\
|S'_1(N)|&|S'_2(N)|&\cdots&|S'_m(N)|
\end{array}\right].$$
{\bf Step 4.} Find usual product of the matrices.\\
{\bf Step 5.} Find complex fuzzy soft max-min decision making matrix $(CSMmDM).$\\
{\bf Step 6.} Find optimum fuzzy set on U.\\

\subsection*{Application}
Let us assume that the set of four signals $U = \{\nu_1, \nu_2, \nu_3, \nu_4\}.$ Now each of these signals is sampled four times. Let $R$ be the given known reference signal. Each signal is compared with the reference signal in order to get the high degree of resemblance with the reference signal $R.$ First use step $(1),$ then we obtain the matrix $A$ by setting the signals along
column and their four times sampling along row. Similarly we will obtain the matrix $B.$\\
Now on the basis of steps algorithm $(1 - 3)$ defined above we discuss an example.
$${A}=\left(\begin{array}{cccc}
0.1 & 0.0 & 0.3 & 0.3\\
0.3 & 0.0 & 0.2 & 0.1\\
0.2 & 0.0 & 0.3 & 0.2\\
0.3 & 0.0 & 0.1 & 0.3
\end{array}\right)$$
$${B}=\left(\begin{array}{cccc}
0.0 & 0.2 & 0.1 &0.4\\
0.0 & 0.1 & 0.3 &0.3\\
0.0 & 0.4 & 0.1 &0.2\\
0.0 & 0.3 & 0.2 &0.1
\end{array}\right).$$
{\bf Step 4.} Now the product of fuzzy soft matrices ${A}$ and ${B}$ is given as:
$${A}\ast{B}=\left(\begin{array}{cccc}
0.00 & 0.23 & 0.01 & 0.13\\
0.00 & 0.17 & 0.07 & 0.17\\
0.00 & 0.22 & 0.09 & 0.18\\
0.00 & 0.19 & 0.10 & 0.18
\end{array}\right).$$
It is simple usual matrix multiplication.\\
{\bf Step 5.} To calculate $Mm\left[{A}\ast{B}\right]= [D_{i1}] = [d_{i1}],$ we have to find $[D_{i1}]$ for all $i\in
\{1, 2, 3, 4\},$ such that
$$d_{i1} = min \{t_{k1}\} = min \{t_{11}, t_{21}, t_{31}, t_{41}\},\; for \;k \in \{1, 2, 3, 4\}.$$
Now let us calculate $[D_{11}] = (d_{11})$ for fixed $j = 1.$
$$d_{11} = min \{t_{k1}\} = min \{t_{11}, t_{21}, t_{31}, t_{41}\}.$$
Here we have to find ${t}_{k1}$ for all $k \in\{1, 2, 3, 4\}$ and $t_{k1}$ for all $k\in \{1, 2, 3, 4\}.$ Let us find
$\overline{t}_{k1}$ for $k\in\{1, 2, 3, 4\}$ as:
$$t_{11 }= 0.00, t_{21} = 0.00, t_{31} = 0.00, t_{41} = 0.00.$$
Now
$d_{11} = min \{0.00, 0.00, 0.00, 0.00\} = 0.00$
So $(D_{11}) = (d_{11}) = (0.00).$ Similarly we can find $(D_{21}) = (d_{21}) , (D_{31}) = (d_{31})$ and $(D_{41}) = (d_{41}),$ where $$d_{21} = min\{t_{12}, t_{22}, t_{32}, t_{42}\} = min \{0.23, 0.17, 0.22, 0.19\} = 0.17$$
So $(D_{21}) = (d_{21}) = 0.17.$
$$d_{31} = min \{0.10, 0.07, 0.09, 0.10\} = 0.07.$$
Therefore $D_{31} = d_{31} = 0.07.$\\
Also $$d_{41} = min \{0.13, 0.17, 0.18, 0.18\} = 0.13$$
So $(D_{41}) = d_{41} = 0.13.$\\
Finally we can obtain fuzzy soft max-min decision fuzzy soft matrix as:
$$mM\left({A}\ast{B}\right)=(D_{i1}) = (d_{i1})$$
$$=\left(\begin{array}{c}
(D_{11}) = (d_{11})\\
(D_{21}) = (d_{21})\\
(D_{31}) = (d_{31})\\
(D_{41}) = (d_{41})
\end{array}\right)=\left(\begin{array}{c}
0.00\\
\textbf{0.17}\\
0.07\\
0.13
\end{array}\right).$$
{\bf Step 6.} Finally we can find an optimum fuzzy set on $[U] = (u)$.
$${opt}Mm(A\ast B) (\nu) = \{\textbf{0.17}/\nu_2, 0.07/\nu_3, 0.13/\nu_4\}.$$
Hence, identify signal $\nu_2$ as the reference signal $R.$
\ \\ \ \\
\section{Decision Making Algorithm in Signal and Systems using Fourier Transform} 
Complex fuzzy soft matrices can represent various characteristics of signals, such as frequency content, amplitude, phase, and noise levels, in a comprehensive and flexible manner. \\
Each element of the matrix can represent a specific signal attribute with fuzzy and soft characteristics, allowing for the incorporation of uncertainty in signal descriptions.

\subsection*{Applications}
\begin{definition}
	The discrete Fourier transform transforms a sequence of $N$ complex numbers $\{x_n(N)\}$ into another sequence of complex numbers $X(k),$ which is defined by
	$$\;\;\;\;\;\;\;\;\;\;\;\;\;\;\;\;X(k)=\sum^{N-1}_{n=0}x_i(n)e^{i(-{2\pi\over N}kn)}\;\;k,n\in\{0,1,2,\cdots,N-1\}.\;\;\;\;\;\;\;\;\;\;\;\;\;\;\;\;\;\;\;\;\;(1)$$
	The inverse discrete Fourier transform of $(1)$ is defined as:
	$$\;\;\;\;\;\;\;\;\;\;\;\;\;\;\;\;x_i(n)={1\over N}\sum^{N-1}_{k=0}X[k]e^{i({2\pi\over N}kn)}\;\;k,n\in\{0,1,2,\cdots,N-1\}.\;\;\;\;\;\;\;\;\;\;\;\;\;\;\;\;\;\;\;(2)$$ which is also $N$-periodic and $X_k$ has different values.
\end{definition}
\begin{remark}
	We take a particular case, that is, $X(k)$ is restricted to a closed interval $[0, 1]$ because, in complex fuzzy set, the amplitude term
	has all the values in the closed interval $[0,1].$
\end{remark}
In the following, an algorithm is develop in signals and systems for the identification of a reference signal received by a particular receiver using complex fuzzy soft sets.\\

\subsection*{Algorithm}
Let $n$ be different electromagnetic signals, $x_1(n), x_2(n), x_3(n),\cdots, x_n(n)$ which have been received by a particular receiver.
Each of these signals is noted $N$ different times. Let $x_i(n)$ be the $i^{th}\; (1 \leq n \leq N)$ signal. The inverse discrete Fourier transform of this $i^{th}$ signal is given as
$$\;\;\;\;\;\;\;\;\;\;\;\;\;\;\;\;x_i(n)={1\over N}\sum^{N-1}_{k=0}X[k]e^{i({2\pi\over N}kn)}\;\;k,n\in\{0,1,2,\cdots,N-1\}.\;\;\;\;\;\;\;\;\;\;\;\;\;\;\;\;\;\;\;(3)$$
We restrict the range of $X[k]$ as $0\leq X[k] \leq1 \;(0 \leq k \leq N-1).$ Here, $X[k]$ is known as amplitude term and
${2\pi\over N}kn=\omega_s(q)$ is known as phase term and the first one having the range as real numbers.\\
Thus, a general signal representing by equations $(3)$ is model for signal representation using a complex fuzzy set.\\
We use the complex fuzzy set in signals and systems utilizing a new kind of matrix to identifies a reference signal out of large signals detected by a digital receiver. For this, we have a reference signal $r$. This reference signal $r$ is noted $N$ times.
The IDFT( inverse discrete Fourier transform) of this reference signal $r$ is
$$\;\;\;\;\;\;\;\;\;\;\;\;\;\;\;\;r(n)={1\over N}\sum^{N-1}_{k=0}X'[k]e^{i({2\pi\over N}kn)}\;\;k,n\in\{0,1,2,\cdots,N-1\}.\;\;\;\;\;\;\;\;\;\;\;\;\;\;\;\;\;\;\;$$
where $X'[k]\in [0,1];\;(0\leq k\leq N-1)$. \\
To compare the similarities between two signals, we apply the following methods.\\
\begin{itemize}
	\item \textbf{Step 1} \\
	Expand $$\;\;\;\;\;\;\;\;\;\;\;\;\;\;\;\;x_i(n)={1\over N}\sum^{N-1}_{k=0}X[k]e^{i({2\pi\over N}kn)}\;\;k = \{0,1,2,\cdots,N-1\}.\;\;\;\;\;\;\;\;\;\;\;\;\;\;\;\;\;\;\;$$ we get
	
	\[
	\begin{aligned}
		x_i(n) &= X[0] e^{i \left( -\frac{2\pi}{N} n(0) \right)} + X[1] e^{i \left( -\frac{2\pi}{N} n(1) \right)} + X[2] e^{i \left( -\frac{2\pi}{N} n(2) \right)} + \ldots + X[N-1] e^{i \left( -\frac{2\pi}{N} n(N-1) \right)} \\
		&= X[0] e^{i \left( -\frac{2\pi}{N} n(0) \right)} + X[1] e^{i \left( -\frac{2\pi}{N} n(1) \right)} + X[2] e^{i \left( -\frac{2\pi}{N} n(2) \right)} + \ldots + X[N-1] e^{i \left( -\frac{2\pi}{N} n(N-1) \right)}
	\end{aligned}
	\]
	
	\[
	x_i(n) = \frac{1}{N} \left[ X[0] \cdot 1 + X[1] e^{i \left( \frac{2\pi}{N} n(1) \right)} + X[2] e^{i \left( \frac{2\pi}{N} n(2) \right)} + \ldots + X[N-1] e^{i \left( \frac{2\pi}{N} n(N-1) \right)} \right] \qquad (5)
	\]
	From equation (5), we get $N$ samples by putting $n = 0, 1, 2, 3, \ldots, N - 1$.\\
	For $n = 0$, we have
	\[ x_i(0) = \frac{1}{N} \left[ X[0] \cdot 1 + X[1] e^{i \frac{2\pi}{N} \cdot (0)(1)} + X[2] e^{i \frac{2\pi}{N} \cdot (0)(2)} + \ldots + X[N-1] e^{i \frac{2\pi}{N} \cdot (0)(N-1)} \right] \]
	
	\[ = \frac{1}{N} \left[ X[0] \cdot 1 + X[1] \cdot 1 + X[2] \cdot 1 + \ldots + X[N-1] \cdot 1 \right] \qquad (6)\] 
	For $n = 1$, we have
	
	\[ x_i(1) = \frac{1}{N} \left[ X[0] \cdot 1 + X[1] e^{i \frac{2\pi}{N}(1)(1)} + X[2] e^{i \frac{2\pi}{N}(1)(2)} + \ldots + X[N - 1] e^{i \frac{2\pi}{N}(1)(N - 1)} \right] \]
	
	\[ = \frac{1}{N} \left[ X[0] \cdot 1 + X[1] e^{i 2\pi (1)} + X[2] e^{i 2\pi (2)} + \ldots + X[N - 1] e^{i 2\pi (N - 1)} \right]. \qquad (7) \]
	For $n = 2$, we have
	
	\[ x_i(2) = \frac{1}{N} \left[ X[0] \cdot 1 + X[1] e^{i \frac{2\pi}{2} (2)(1)} + X[2] e^{i \frac{2\pi}{2} (2)(2)} + \ldots + X[N - 1] e^{i \frac{2\pi}{2} (2)(N - 1)} \right] \]
	
	\[ = \frac{1}{N} \left[ X[0] \cdot 1 + X[1] e^{i 2\pi (2)} + X[2] e^{i 2\pi (4)} + \ldots + X[N - 1] e^{i 2\pi \cdot 2(N - 1)} \right]. \qquad (8) \]
	
	Continuing this process, for $n = N - 1$, we have
	\[ x_i(N - 1) = \frac{1}{N} \left[ X[0] \cdot 1 + X[1] e^{i \frac{2\pi}{2} (N-1)(1)} + X[2] e^{i \frac{2\pi}{2} (N-1)(2)} + \ldots + X[N - 1] e^{i \frac{2\pi}{2} (N-1)(N-1)} \right] \]
	
	\[ = \frac{1}{N} \left[ X[0] \cdot 1 + X[1] e^{i \frac{2\pi}{2} (N-1)} + X[2] e^{i \frac{2\pi}{2} 2(N-1)} + \ldots + X[N - 1] e^{i \frac{2\pi}{2} (N-1)^2} \right]. \qquad (9) \]
	A similar argument repeats for the reference signal $r(n)$, where we obtain $N$ samples of the reference signal by putting $n = 0, 1, 2, 3, \ldots, N - 1$. \\
	\item \textbf{Step 2}\\
	Now, we find the cross product of these $N$ samples of the signal $x_i(n)$ and the reference signal $r(n)$, where $n = 0, 1, 2, \ldots, N - 1$.
	\[
	\begin{aligned}
		x_i(0) \ast r(0) &= \frac{1}{N^2} \Bigg[\dfrac{\min\{X[0] \times X'[0]\}e^{{i} \min\left\{ \frac{2\pi}{N} N(0), \frac{2\pi}{N} N(0) \right\}}}{(0, 0)} +\\ &\dfrac{\min\{X[0] \times X'[0]\}e^{{i} \min\left\{ \frac{2\pi}{N} N(0), \frac{2\pi}{N} N(0)(1) \right\}}}{(0, 1)} + \ldots+\\&\dfrac{\min\{X[0] \times X'[N-1]\}e^{{i} \min\left\{ \frac{2\pi}{N} N(0), \frac{2\pi}{N} N(0)(N-1) \right\}}}{(0, N-1)} +\\&\dfrac{\min\{X[1] \times X'[0]\}e^{{i} \min\left\{ \frac{2\pi}{N} N(0)(1), \frac{2\pi}{N} N(0)(0) \right\}}}{(1, 0)} + \\& \dfrac{\min\{X[1] \times X'[1]\}e^{{i} \min\left\{ \frac{2\pi}{N} N(0)(1), \frac{2\pi}{N} N(0)(0) \right\}}}{(1, 1)}  + \ldots+\\& \dfrac{\min\{X[1] \times X'[N-1]\}e^{{i} \min\left\{ \frac{2\pi}{N} (0)(1), \frac{2\pi}{N} N(0)(N-1) \right\}}}{(1, N-1)} + \ldots +\\&\dfrac{\min\{X[N-1] \times X'[0]\}e^{{i} \min\left\{ \frac{2\pi}{N} (0)(N-1), \frac{2\pi}{N} (0)(0) \right\}}}{(N-1, 0)} + \ldots +\\&\dfrac{\min\{X[N-1] \times X'[N-1]\}e^{{i} \min\left\{ \frac{2\pi}{N} (0)(N-1), \frac{2\pi}{N} (0)(N-1) \right\}}}{(0,1)}\Bigg]
	\end{aligned}
	\]
	Compute the absolute value of every term of $x_i(0) \times r(0)$ and divide by $N^2$. Now, we take the maximum value from the absolute values of $x_i(0) \times r(0)$. A similar process repeats for $x_i(1) \times r(1)$, $x_i(2) \times r(2)$, ..., $x_i(N-1) \times r(N-1)$. Now, we develop the column matrix from all these max values of $x_i(0) \times r(0)$, $x_i(1) \times r(1)$, $x_i(2) \times r(2)$, ..., $x_i(N-1) \times r(N-1)$; that is,
	
	\[
	\begin{bmatrix}
		\max |x_i(0) \times r(0)| \\
		\max |x_i(1) \times r(1)| \\
		\max |x_i(2) \times r(2)| \\
		\vdots \\
		\max |x_i(N - 1) \times r(N - 1)| \\
	\end{bmatrix}    \hspace{3cm} [10]
	\]
	Similarly, we develop the column matrix from the cross product of the $N - 1$ samples of the signal $x_j(n)$ and the reference signal $r(n)$, that is, we have
	\[
	\begin{bmatrix}
		\max |x_j(0) \times r(0)| \\
		\max |x_j(1) \times r(1)| \\
		\max |x_j(2) \times r(2)| \\
		\vdots \\
		\max |x_j(N - 1) \times r(N - 1)| \\
	\end{bmatrix}  \hspace{3cm} [11]
	\]
	
	\item \textbf{Step 3}\\ For identification of the reference signal of $x_i(n)$ and $x_j(n)$; $n = 0, 1, 2, \ldots, N - 1$, we again take the maximum value of the column matrix (10) and (11). If the max value of column matrix (10) is greater than the max value of the column matrix (11), then $x_i(n)$ shows the reference signal. If the max value of the column matrix (11) is greater than the max value of the column matrix (10), then $x_j(n)$ is the reference signal.\\
%	\vspace{2cm}
\par Figure 1 below shows the graphs of the N-sampled reference signal, \(i-th \) signal, and \(j-th\) signal. Note that these are
discrete signal graphs.\\
\end{itemize}
\begin{figure}[h!]
	\centering
	\includegraphics[width=16cm]{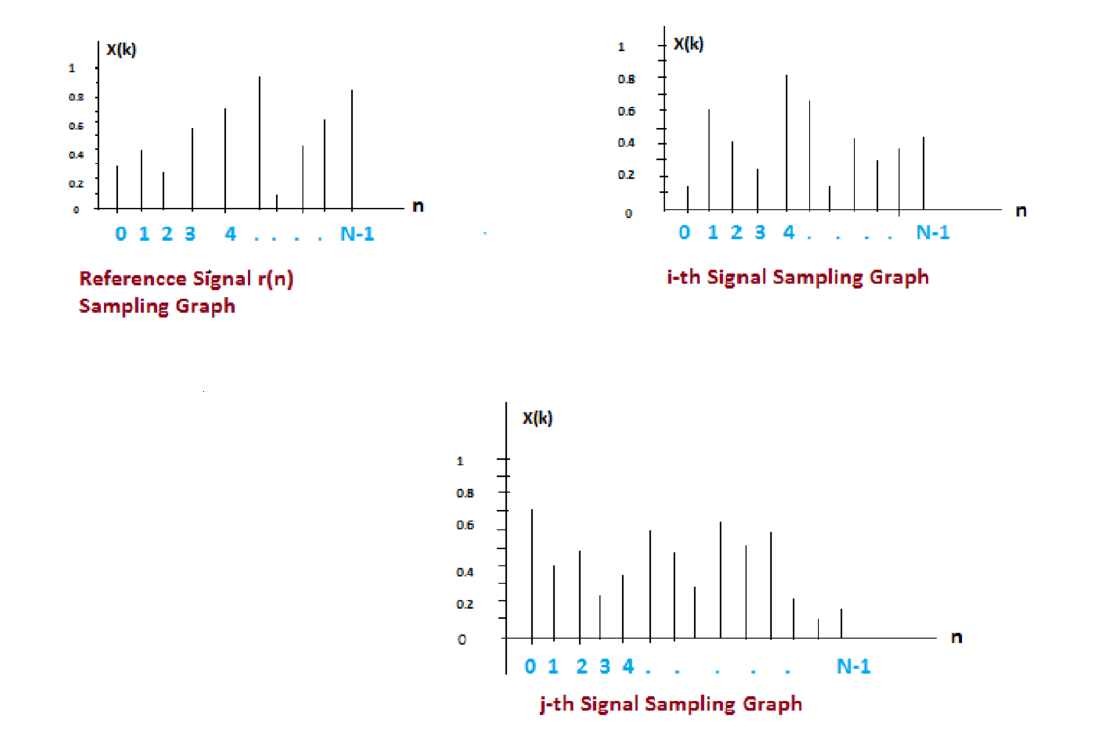}
	\caption{}
	\label{fig:11}
\end{figure}
\begin{figure}[h!]
	\centering
	\includegraphics[width=9cm]{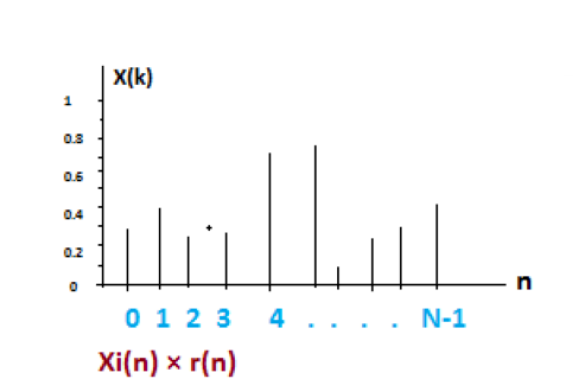}
	\caption{Graph of \(i-th\) signal x reference signal, "x" denotes cross product of CFSs}
	\label{fig:12}
\end{figure}
\begin{figure}[h!]
	\centering
	\includegraphics[width=9cm]{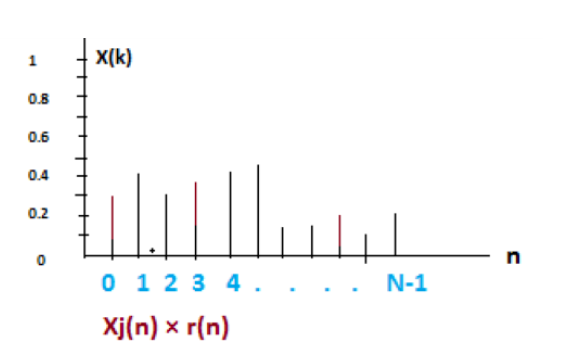}
	\caption{Graph of \(j-th\) signal x reference signal, "x" denotes cross product of CFSs}
	\label{fig:13}
\end{figure}
%$\vspace {3cm}$
\par From Figure 2 and Figure 3, it is clear that the $i$-th signal shows a high degree of resemblance to the reference signal. Thus the $i$-th signal is the reference signal \\
Note that our proposed model discussed the resemblance of signals from different sources. If the electromagnetic signals are from one source, then there is no need to check the resemblance.\\
We may apply different techniques and operations to solve noise signal problems. If the Fourier transform of noise signals is in discrete form, then we can apply our proposed algorithm. If the Fourier transform of noise signals is in continuous form, then we may apply continuous Fourier transform instead of discrete Fourier transform.

\begin{example}
	Assume that two different electromagnetic signals, $x_1(n)$, $x_2(n)$ have been received by a receiver. Each of these signals is sampled two times. Let $r(n)$ be the reference signal. The discrete Fourier transform of the signal $x_i(n)$; $n = 0, 1$ and reference signal $r(n)$ for $N = 2$ is
\end{example}
\[
x_i(n) = \frac{1}{2} \sum_{k=0}^{1} X[k]e^{\frac{2\pi}{2}kn}; \quad k, n = 0, 1, \quad\quad (1)
\]

where $X[k] \in [0, 1]$. Also,

\[
r(n) = \frac{1}{2} \sum_{k=0}^{1} X'[k]e^{\frac{2\pi}{2}nk}; \quad n, k = 0, 1, \quad\quad (2)
\]
where $X'[k] \in [0, 1]$. For $k = 0, 1$, equation (1) becomes
\[
x_i(n) = \frac{1}{2} [X[0]e^{i\frac{2\pi}{2}n(0)} + X[1]e^{i\frac{2\pi}{2}n(1)}] = \frac{1}{2} [X[0]\cdot 1 + X[1]\cdot e^{i\pi n}]. \quad\quad (3)
\]
Now put $n = 0$ and $i = 1$ in (3), we have
\[
x_1(0) = \frac{1}{2} [X[0]e^{i\pi(0)(0)} + X[1]e^{i\pi(0)(1)}] = \frac{1}{2} [X[0]e^{i(0)} + X[1]e^{i(0)}]. \quad \quad(3.1)
\]
Put $n = 1$, we have
\[
x_1(1) = \frac{1}{2} [X[0]e^{i\pi(0)(1)} + X[1]e^{i\pi(1)}] = \frac{1}{2} [X[0]e^{i\pi(0)} + X[1]e^{i\pi(1)}]. \quad\quad (3.2)
\]
A similar process for the reference signal, we have
\[
r(0) = \frac{1}{2} [X'[0]e^{i\pi(0)(0)} + X'[1]e^{i\pi(0)(1)}] = \frac{1}{2} [X'[0]e^{i\pi(0)} + X'[1]e^{i\pi(0)}]. \quad\quad (2.1)
\]
\[
r(1) = \frac{1}{2} [X'[0]e^{i\pi(0)(1)} + X'[1]e^{i\pi(1)(1)}] = \frac{1}{2} [X'[0]e^{i\pi(0)} + X'[1]e^{i\pi(1)}] \quad \quad(2.2)
\]
We take the amplitude term of the signal $x_1(0)$, $x_1(1)$, $r(0)$, and $r(1)$ respectively as
\[
X[k] = \begin{cases}
	0.7 & \text{; } k = 0 \\
	0.4 & \text{; } k = 1 ,\\
\end{cases} \]
\[
X[k] = \begin{cases}
	1.0 & \text{; } k = 0 \\
	0.3 & \text{; } k = 1, \\
\end{cases}\]
\[X'[k] = \begin{cases}
	0.1 & \text{; } k = 0 \\
	0.9 & \text{; } k = 1, \\
\end{cases}
\]
\[
X'[k] = \begin{cases}
	0.6 & \text{; } k = 0 \\
	0.5 & \text{; } k = 1. \\
\end{cases}
\]
Therefore, (2.3), (2.4), (3.3), and (3.4) become respectively
\[
r(0) = \frac{1}{2}[0.1\cdot e^{i\pi(0)(0)} + 0.9\cdot e^{i\pi(0)(1)}], \quad (2.3)
\]
\[
r(1) = \frac{1}{2}[0.6\cdot e^{i\pi(0)(1)} + 0.5\cdot e^{i\pi(1)(1)}]. \quad (2.4)
\]
\[
x_1(0) = \frac{1}{2}[0.7\cdot e^{i\pi(0)(0)} + 0.4\cdot e^{i\pi(0)(1)}]. \quad (3.3)
\]
\[
x_1(1) = \frac{1}{2}[1.0\cdot e^{i\pi(0)(1)} + 0.3\cdot e^{i\pi(1)}]. \quad (3.4)
\]
Now, we find the cross product of $x_1(n)$ and the reference signal $r(N)$; that is
\[
\begin{aligned}
	x_1(0) \times r(0) &= \frac{1}{4} [ \min\{0.1, 0.7\}e^{\text{i}  min\{0, 0\}} + \min\{0,0.5\}e^{\text{i} min\{0, 0\}} + \\& \min\{1.1,0.4\}e^{\text{i} min\{0, 0\}} + \min\{0.9,0.7\}e^{\text{i} min\{0.9, 0.4\}}] \\	&= \frac{1}{4} [0.1(1) + 0.1(1) + 0.7(1) + 0.4(1)]\\
	&= \frac{1}{4} [0.1 + 0.13 + 0.7 + 0.4].
\end{aligned}
\]
Now,
\begin{align*}
	\max |x_1(0) \times r(0)| &= \max \left[ |\frac{0.1}{4}| + |\frac{0.1}{4}| + |\frac{0.7}{4}| + |\frac{0.4}{4}| \right] \\&= \max[0.03 + 0.03 + 0.18 + 0.10] = 0.18.
\end{align*}
Thus, \[max | x_1 (0) \times\ r(0)|= 0.18\]
Also,
\begin{align*}
	x_1(1) \times r(1) &= \frac{1}{4} [\min\{0.6,1.0\}e^{\text{i} min\{0, 0\}} + \min\{0.6,0.3\}e^{\text{i} min\{0, \pi\}} + \\ &\min\{0.5,1.0\}e^{\text{i} min\{\pi, 0\}} + \min\{0.5,0.3\}e^{\text{i} min\{\pi , \pi\}} ]\\&
	=\frac{1}{4} [0.6(1) + 0.3(1) + 0.5(1) + 0.3(-1)]\\&
	=\frac{1}{4} [0.6 + 0.3 + 0.5 - 0.3] = \frac{0.6 + 0.3 + 0.3 - 0.3}{4}
\end{align*}
Now,
\begin{align*}
	\max |x_1(1) \times r(1)| &= \max \left[ |\frac{0.6}{4}| + |\frac{0.3}{4}| + |\frac{0.5}{4}| + |\frac{-0.3}{4}| \right] \\&= \max[0.15 + 0.08 + 0.13 + 0.08] = 0.15.
\end{align*}
Thus, \[max | x_1 (1) \times\ r(1)|= 0.15\]
Now, we develop the column matrix from all these max values of $x_1(0) \times r(0)$ and $x_1(1) \times r(1)$,\\ $x_n(2) \times r(2)$; that is
\[
\begin{bmatrix}
	\max |x_1(0) \times r(0)| \\
	\max |x_1(1) \times r(1)| \\
\end{bmatrix}
=
\begin{bmatrix}
	0.18 \\
	0.15 \\
\end{bmatrix}
\hspace{3cm} (A)
\]
Now consider the signal
\[
x_2(n) = \frac{1}{2} \sum_{k=0}^{1} X[k]e^{i\frac{2\pi}{2}kn}; \quad k, n = 0, 1, \quad (4)
\]
where $X[k] \in [0, 1]$. For $k = 0, 1$, equation (1) becomes
\[
x_2(n) = \frac{1}{2} [X[0]e^{i\frac{2\pi}{2}n(0)} + X[1]e^{i\frac{2\pi}{2}n(1)}] = \frac{1}{2} [X[0]\cdot 1 + X[1]\cdot e^{i\pi n}]. \quad (4.1)
\]
Now put $n = 0$ and $1$ in (4.1), we get
\[
x_2(0) = \frac{1}{2} [X[0]e^{i\pi(0)(0)} + X[1]e^{i\pi(0)(1)}], \quad (4.2)
\]
\[
x_2(1) = \frac{1}{2} [X[0]e^{i\pi(0)(1)} + X[1]e^{i\pi(1)}]. \quad (4.3)
\]
We take the amplitude term of the signal $x_2(0)$, $x_2(1)$, respectively as:
\[
X[k] = \begin{cases}
	0.9 & \text{; } k = 0 \\
	0.6 & \text{; } k = 1, \\
\end{cases} \]
\[
X[k] = \begin{cases}
	0.8 & \text{; } k = 0 \\
	0.1 & \text{; } k = 1 \\
\end{cases}
.
\]
Therefore, (4.2) and (4.3) become as:
\[
x_2(0) = \frac{1}{2} [0.9e^{i\pi(0)(0)} + 0.6e^{i\pi(0)(1)}], \quad (4.4)
\]
\[
x_2(1) = \frac{1}{2} [0.8e^{i\pi(0)(1)} + 0.1e^{i\pi(1)}]. \quad (4.5)
\]
Now, we find the cross product of $x_2(n)$ and $r(n)$, that is,
\[
\begin{aligned}
	x_2(0) \times r(0) &= \frac{1}{4} [ \min\{0.1,0.9\}e^{\text{i} min\{0, 0\}} + \min\{0.1,0.6\}e^{\text{i} min\{0, 0\}} \\&+ \min\{0.9,0.9\}e^{\text{i} min\{0, 0\}} + \min\{0.9,0.6\}e^{\text{i} min\{0, 0\}} ] \\&= \frac{1}{4} [0.1(1) + 0.1(1) + 0.9(1) + 0.6(1)]  \\
	&= \frac{1}{4} [0.1 + 0.1 + 0.9 + 0.6] = \frac{0.1 + 0.1 + 0.9 + 0.6}{4}.
\end{aligned}
\]
Now,
\begin{align*}
	\max |x_2(0) \times r(0)| &= \max \left[ \frac{|0.1|}{4} + \frac{|0.1|}{4} + \frac{|0.9|}{4} + \frac{|0.6|}{4} \right] \\&= \max[0.03 + 0.03 + 0.23 + 0.15] = 0.23.
\end{align*}Thus,
\[max | x_2 (0) \times\ r(0)|= 0.23\]
Now,
\begin{align*}
	x_2(1) \times r(1) &= \frac{1}{4} [ \min\{0.6,0.8\}e^{\text{i} min\{0, 0\}} + \min\{0.6,0.1\}e^{\text{i} min\{0, \pi\}} +\\& \min\{0.5,0.8\}e^{\text{i} min\{\pi, 0\}} + \min\{0.5,0.1\}e^{\text{i} min\{\pi , \pi\}} ] \\&=\frac{1}{4} [0.6(1) + 0.1(1) + 0.5(1) + 0.1(-1)]\\&=\frac{1}{4} [0.6 + 0.1 + 0.5 - 0.1] = \frac{0.6 + 0.1 + 0.5 - 0.1}{4}
\end{align*}
Now,
\begin{align*}
	\max |x_2(1) \times r(1)| &= \max [ |\frac{0.6}{4}| + |\frac{0.1}{4}| + |\frac{0.5}{4}| + |\frac{-0.1}{4}|] \\&= \max[0.15 + 0.03 + 0.13 + 0.03] = 0.15
\end{align*} 
Thus,
\[max | x_2 (1) \times\ r(1)|= 0.15\]
Now, we develop the column matrix from all these max values of $x_2(0) \times r(0)$ and\\ $x_2(1) \times r(1)$, that is,
\[
\begin{bmatrix}
	\max |x_2(0) \times r(0)| \\
	\max |x_2(1) \times r(1)| \\
\end{bmatrix}
=
\begin{bmatrix}
	\textbf{0.23} \\
	\textbf{0.15} \\
\end{bmatrix} \hspace{3cm} (B)
\]
From (A) and (B), we have
\[ max \begin{bmatrix}
	\max |x_1(0) \times r(0)| \\
	\max |x_1(1) \times r(1)| \\
\end{bmatrix} > max \begin{bmatrix}
	\max |x_2(0) \times r(0)| \\
	\max |x_2(1) \times r(1)| \\
\end{bmatrix}\]
Thus, the signal $x_2(n)$ which is matrix B identifies the reference signal.

\ \\ \ \\

\section*{Conclusion}
One of the key problems in signal processing is how to select a suitable model. 
In this paper, we have discussed the application of complex fuzzy soft sets in signals and systems via the cross product of complex fuzzy soft matrices and Fourier transform. In this application,  an algorithm for the identification of a reference signal out of large interest signals detected by a digital receiver was presented. We utilised the cross product of the known signal and unknown signals and take the maximum value. We compared it with all the cross products of the known and unknown signals. We again take the maximum value among all of these values and observed that the unknown signal shows a high degree of resemblance to the reference signal. Thus we identified one unknown signal as a reference signal $R$ among the several signals detected by the receiver. It was noted that, the Fourier transform is better because it gave a higher optimal value and as a result, there was a better reference signal $R.$  
\ \\ \ \\ \ \\

%%%%%%%%%%%%%%%%%%%%%%%%%%%%%%%%%%%%%%%%%%%%%%%%%%%%%%%%%%%%%%

\end{document}